\documentclass[twoside]{kybernetika}
\usepackage{graphicx, amssymb, amsmath}
\usepackage{hyperref}

\newtheorem{theorem}{Theorem}[section]
\newtheorem{lemma}[theorem]{Lemma}

\newtheorem{example}[theorem]{Example}

\def\msn{\medskip\noindent}

\def\RR{\overline{\mathbb{R}}}

\def\R{\mathbb{R}}

\def\spec{\sigma}

\def\supp{\operatorname{supp}}
\def\spann{\operatorname{span}_{\oplus}}

\begin{document}
\pagestyle{myheadings}

\title{On the problem $Ax=\lambda Bx$ in max algebra:\\ every
system of intervals is a spectrum}

\author{Serge{\u{\i}} Sergeev}

\contact{Serge{\u{\i}}}{Sergeev}{CMAP, \'Ecole Polytechnique, 91128 Palaiseau C\'edex, France} {sergeev@cmap.polytechnique.fr}

\markboth{S. Sergeev} {A note on the two-sided eigenproblem}

\maketitle

\begin{abstract}
We consider the two-sided eigenproblem $A\otimes x=\lambda\otimes B\otimes x$
over max algebra. It is shown that any finite system of 
real intervals and points can be represented as spectrum of this eigenproblem.
\end{abstract}

\keywords{Extremal algebra, tropical algebra, generalized eigenproblem}

\classification{15A80,15A22, 91A46,93C65}

\section{INTRODUCTION}

Max algebra is the analogue of linear algebra
developed over the max-plus semiring, which is
the set $\RR=\R\cup\{-\infty\}$ equipped with the operations
of ``addition'' $a\oplus b:=\max(a,b)$ and 
``multiplication'' $a\otimes b:=a+b$. 
This basic arithmetics is naturally extended to matrices and vectors.
In particular, for matrices $A=(a_{ij})\in\RR^{n\times m}$ and 
$B=(b_{ij})\in\RR^{m\times k}$, their ``product''
$A\otimes B$ is defined by the rule 
$(A\otimes B)_{ij}=\bigoplus_{l=1}^m a_{il}\otimes b_{lj}$, for all $i=1,\ldots,n$
and $j=1,\ldots,k$.

One of the best studied problems in max algebra is the ``eigenproblem'':
for given $A\in\RR^{n\times n}$ find $\lambda\in\RR$ and $x\in\RR^n$ with
at least one finite entry, such that $A\otimes x=\lambda\otimes x$.
This problem is very important for max-algebra and its applications 
\cite{ABG-06,BCOQ,But-03,But:10, CG:79, HOW:05}. The theory of this problem has much in common
with its counterpart in the nonnegative matrix algebra. In particular,
there is exactly one eigenvalue (``max-algebraic Perron root'') in the irreducible
case, and in general, there may be several eigenvalues which correspond
to diagonal blocks of the Frobenius normal form. There are efficient
algorithms for computing both eigenvalues and eigenvectors \cite{HOW:05,Coc-98,ED-01}.

We will consider the following generalization of the max algebraic
eigenproblem: 
\begin{equation}
\label{twosided-eig}
A\otimes x=\lambda\otimes B\otimes x,
\end{equation}
where $A,B\in\RR^{n\times m}$. The set of $\lambda\in\RR$ such that there
exists $x$ satisfying \eqref{twosided-eig}, with at least one finite entry,
will be called the {\em spectrum} of \eqref{twosided-eig} and denoted by 
$\spec(A,B)$.

This problem is of interest as an analogue of matrix pencils in nonnegative matrix algebra,
as studied in McDonald et al.~\cite{McD+98}, Mehrmann et al.~\cite{MNV-08}. Note that matrix pencils in linear algebra are very well-known, see Gantmacher~\cite{Gan:59} for basic reference, and their
applications in control go back to Brunovsky~\cite{Bru-70}.  

Problem~\eqref{twosided-eig} can also be considered as a parametric extension
of two-sided systems $A\otimes x=B\otimes x$. Importantly, such systems 
can be solved algorithmically, see Cuninghame Green and Butkovi\v{c}~\cite{CGB-03}.

Unlike the eigenproblem $A\otimes x=\lambda\otimes x$, 
the two-sided version does not 
seem to be well-known. Some results have been obtained by Binding and Volkmer~
\cite{BV-07}, and Cuninghame-Green and Butkovi{\v{c}}\cite{CGB-08}, mostly for special cases when both matrices are square,
or when $A=B\otimes Q$. See also Butkovi\v{c}~\cite{But:10}. In the latter case, it may be possible
to reduce \eqref{twosided-eig} to $Q\otimes x=\lambda\otimes x$. In general,
however, it is nontrivial to decide whether the spectrum is nonempty, and 
some particular conditions have been studied by topological methods~\cite{BV-07}. 

Further, the spectrum of \eqref{twosided-eig} can be much richer, it may include 
intervals. Gaubert and Sergeev~\cite{GS-10} came up with a general approach to the problem
representing it in terms of parametric min-max functions and mean-payoff games, which allows to 
identify the whole spectrum in pseudo-polynomial time.
The purpose of this note is more modest, it is to provide an example showing that 
any system of intervals and points can be 
realized as the spectrum of \eqref{twosided-eig}. 

Let us note a possible application of \eqref{twosided-eig} in scheduling in the spirit of Cuninghame-Green~\cite{CG:79}. See also Burns~\cite{Bur:91}. 
Suppose that the products $P_1,\ldots,P_n$ are prepared using $m$ machines (or, say, processors),
where every machine contributes to the completion 
of each product by producing a partial product.
Let $a_{ij}$ be the duration of the work of the $j$th machine needed to 
complete the partial product for $P_i$. Let us denote by $x_j$ the
starting time of the $j$th machine, then all partial products for $P_i$ will
be ready by the time $\max(x_1+a_{i1},\ldots,x_m+a_{im})$. Now suppose that
$m$ other machines prepare partial products for products $Q_1,\ldots,Q_n$,
and the duration and starting times are $b_{ij}$ and $y_j$ respectively.
If the machines are linked then it may be required that $y_j-x_j$ is a constant
time $\lambda$. Now consider a synchronization problem: to find $\lambda$
and starting times of all $2m$ machines so that each pair 
$P_i,Q_i$ is completed at the same time. Algebraically, we have
to solve
\begin{equation}
\label{twosided-detail}
\begin{split}
&\max(x_1+a_{i1},\ldots,x_m+a_{im})=\max(\lambda+x_1+b_{i1},\ldots,\lambda+x_m+b_{im}),\\
&\forall i=1,\ldots,n,
\end{split}
\end{equation}
which is clearly the same as \eqref{twosided-eig}.

\section{PRELIMINARIES} 

We begin with some definitions and notation.
The {\em max algebraic column span} of $A=(a_{ij})\in\RR^{n\times m}$
is defined by
\begin{equation*}
\spann(A)=\left\{\bigoplus_{i=1}^m \alpha_i A_{\cdot i}\mid \alpha_i\in\RR\right\}.
\end{equation*}
For $y\in\RR^n$ denote $\supp(y)=\{i:\; y_i\neq-\infty\}$, and
for $y,z\in\RR^n$ denote
\begin{equation*}
T(y,z):=\arg\min\{y_i-z_i\mid i\in \supp(y)\cap\supp(z)\}.
\end{equation*}

In max algebra, {\em one-sided} systems $A\otimes x=b$ 
can be easily solved, and the solvability criterion is as follows. By $A_{i\cdot}$ (resp. $A_{\cdot i}$)
we denote the $i$th row (resp. the $i$th column) of $A\in\RR^{n\times m}$.

\begin{theorem}[\cite{But-03}, Theorem 2.2]
\label{Vorobyov}
Let $A\in\RR^{n\times m}$ and $b\in\R^n$. The following statements
are equivalent.
\begin{itemize}
\item[1.] $b\in\spann(A)$.
\item[2.] $A\otimes x=b$ is solvable.
\item[3.] $\bigcup_{i=1}^m T(b,A_{\cdot i})=\{1,\ldots,n\}$.
\end{itemize}
\end{theorem}

The author is not aware of any such criterion for two-sided systems
$A\otimes x=B\otimes x$. However, the following {\em cancellation law}
can be useful in their analysis ($a,b,c,d\in\RR$):
\begin{equation}
\label{canclaw}
\begin{split}
& \text{if $a<c$ then}\\
& a\otimes  x\oplus b=c\otimes x\oplus d\quad\Leftrightarrow\quad b=c\otimes x\oplus d.
\end{split}
\end{equation}  
Consider a particular application of this law. In what follows we write $x<y$ also for two
{\em vectors} $x$ and $y$, if $x_i<y_i$ holds for all their components. 

\begin{lemma}
\label{l:A<B}
Let $A,B\in\RR^{n\times m}$ and let $A_{i\cdot}<B_{i\cdot}$
for some $i$. Then $A\otimes x=B\otimes x$ does not have nontrivial
solution.
\end{lemma}
\begin{Proof}
Applying cancellation \eqref{canclaw}, we obtain that the $i$th equation
of $A\otimes x=B\otimes x$ is equivalent to $B_{i\cdot}\otimes x=-\infty$.
Note that all entries of $B_{i\cdot}$ are finite, hence $x_j=-\infty$ for all $j$.
\end{Proof}

When $A,B$ have finite entries only, Lemma \ref{l:A<B} can be used \cite{CGB-08}
to obtain bounds for the spectrum of \eqref{twosided-eig}:
\begin{equation}
\label{ABbounds}
\spec(A,B)\subseteq
[\max\limits_i\min\limits_j(a_{ij}-b_{ij}),
\min\limits_i\max\limits_j (a_{ij}-b_{ij})].
\end{equation}
The cancellation law also allows to replace the finiteness restriction by requiring
that $a_{ij}$ {\em or} $b_{ij}$ is finite for all $i$ and $j$.

It will be also useful that \eqref{twosided-eig}
is equivalent to the following system with separated variables:

\begin{equation}
\label{CDdef}
\begin{split}
C(\lambda)\otimes x &= D\otimes y,\ \text{where}\\
C(\lambda)=
\begin{pmatrix}
A\\
\lambda\otimes B
\end{pmatrix}, &\quad
D=
\begin{pmatrix}
I\\
I
\end{pmatrix},
\end{split}
\end{equation}
and $I=(\delta_{ij})\in\RR^{n\times n}$ denotes the max-plus identity matrix
with entries
\begin{equation}
\delta_{ij}=
\begin{cases}
0, & \text{if $i=j$},\\
-\infty, & \text{if $i\ne j$}.
\end{cases}
\end{equation}

The finite vectors belonging to $\spann(D)$ can be easily described.

\begin{lemma}
\label{spannd}
$z\in\R^{2n}$ belongs to $\spann(D)$ if and only if $z_i=z_{n+i}$ for
all $i=1,\ldots,n$.
\end{lemma}

\section{MAIN RESULTS}

Let $\{[a_i,c_i],\ i=1,\ldots,m\}$ be a 
finite system of intervals on the real line, where
$a_i\leq c_i<a_{i+1}$ for $i=1,\ldots, m-1$, with possibility that $a_i=c_i$.
Define matrices $A\in\R^{2\times 3m}$, $B\in\R^{2\times 3m}$:
\begin{equation}
\label{ABdef}
\begin{split}
A&=
\begin{pmatrix}
\ldots & a_i & b_i & c_i &\ldots\\
\ldots & 2a_i & 2b_i & 2c_i & \ldots
\end{pmatrix},\\
B&=
\begin{pmatrix}
\ldots & 0 & 0 & 0 &\ldots\\
\ldots & a_i & c_i & b_i &\ldots
\end{pmatrix},
\end{split}
\end{equation}
where $b_i:=\frac{a_i+c_i}{2}$. 

\begin{theorem}
\label{mainres-spec}
With $A,B$ defined by \eqref{ABdef},
\begin{equation}
\label{e:mainres-spec}
\spec(A,B)=\bigcup_{i=1}^m [a_i,c_i].
\end{equation}
\end{theorem}
\begin{Proof}
First we show that any $\lambda$ outside the system of intervals is not
an eigenvalue.

\msn {\em Case 1.} $\lambda<a_1$, resp. $\lambda>c_m$. In these cases
$\lambda\otimes B_{1\cdot}<A_{1\cdot}$, resp. 
$\lambda\otimes B_{1\cdot}>A_{1\cdot}$, hence by Lemma \ref{l:A<B} 
$A\otimes x=\lambda\otimes B\otimes x$ cannot hold with nontrivial $x$.

\msn {\em Case 2.} $c_k<\lambda<a_{k+1}$. 
Using cancellation law \eqref{canclaw}, we obtain that the first equation of
$A\otimes x=\lambda\otimes B\otimes x$ is equivalent to
\begin{equation}
\label{1steq}
\bigoplus_{i=k}^{m-1} (a_{i+1}\otimes x_{3i+1}\oplus b_{i+1}\otimes x_{3i+2}\oplus
c_{i+1}\otimes x_{3i+3})=\lambda\otimes\bigoplus_{i=1}^{3k} x_i.
\end{equation}
For the second equation of $A\otimes x=\lambda\otimes B\otimes x$, observe that $2a_i>\lambda+a_i$, $2b_i>\lambda+c_i$
and $2c_i>\lambda+b_i$ for all $i\geq k+1$. After cancellation \eqref{canclaw},
the l.h.s. and the r.h.s. of this equation turn into max-linear forms $u(x)$ and 
$v(x)$ respectively, such that
\begin{equation}
\label{2ndeq}
\begin{split}
u(x)&=v(x),\\
u(x)&\geq\bigoplus_{i=k}^{m-1} (2a_{i+1}\otimes x_{3i+1}\oplus 2b_{i+1}\otimes x_{3i+2}\oplus 2c_{i+1}\otimes x_{3i+3}),\\
v(x)&\leq\lambda\otimes\bigoplus_{i=0}^{k-1} (a_{i+1}\otimes x_{3i+1}\oplus c_{i+1}\otimes x_{3i+2}\oplus b_{i+1}\otimes x_{3i+3}).
\end{split}
\end{equation}
We claim that \eqref{1steq} and \eqref{2ndeq} cannot hold at the same time
with a nontrivial $x$. 
Using that $\lambda<a_{i+1}\leq b_{i+1}\leq c_{i+1}$ for all $i\geq k$, and that the l.h.s. of
\eqref{1steq} attains maximum at a particular term, we deduce from
\eqref{1steq} that
\begin{equation}
\label{inter1}
\bigoplus_{i=k}^{m-1} (2a_{i+1}\otimes x_{3i+1}\oplus 2b_{i+1}\otimes x_{3i+2}\oplus
2c_{i+1}\otimes x_{3i+3})>2\lambda\otimes\bigoplus_{i=1}^{3k} x_i.
\end{equation}
(Note that both sides of~\eqref{1steq} are finite, since all coefficients of $A$ and $B$ are
finite and $x$ is nontrivial.)
The l.h.s. of~\eqref{inter1} is the same as the r.h.s. of the second statement of~\eqref{2ndeq}.
Therefore, combining \eqref{inter1} and the first two statements of \eqref{2ndeq}, we obtain
\begin{equation}
\label{inter2}
v(x)=u(x)>2\lambda\otimes\bigoplus_{i=1}^{3k} x_i,
\end{equation}
Now, since the coefficients $a_{i+1}$, $b_{i+1}$ and $c_{i+1}$ on the r.h.s. of the last statement of \eqref{2ndeq} do not exceed
$\lambda$, the r.h.s. of that statement does not exceed the r.h.s. of 
\eqref{inter2}. But combining \eqref{inter2} with that last statement of \eqref{2ndeq} we obtain
just the opposite. This contradiction shows that 
$A\otimes x=\lambda\otimes B\otimes x$ cannot have nontrivial
solutions in case 2.

Now we prove that any $\lambda$ in the intervals is an eigenvalue, by
guessing a vector that belongs to $\spann(C(\lambda))\cap\spann(D)$.
The columns of $C(\lambda)$ will be denoted by
\begin{equation}
\label{uvwdef}
\begin{array}{c@{{}\quad{}}ccc}
u^i(\lambda)=(a_i & 2a_i & \lambda & a_i+\lambda )^T,\\
v^i(\lambda)=(b_i & 2b_i & \lambda & c_i+\lambda )^T,\\
w^i(\lambda)=(c_i & 2c_i & \lambda & b_i+\lambda )^T.
\end{array}
\end{equation}

\msn {\em Case 3.} $a_i\leq\lambda\leq b_i$. We take
\begin{equation}
\label{ydef1}
\begin{array}{c@{{}\quad{}}ccc}
z^{\lambda}=(0 & \lambda+b_i-a_i & 0 & \lambda+b_i-a_i)^T.
\end{array}
\end{equation}
By Lemma \ref{spannd} $z^{\lambda}\in\spann(D)$. 
It suffices to check that $z^{\lambda}\in\spann(C(\lambda))$.
We write
\begin{align}
\label{e:y/u:1}
T(z^{\lambda},u^i(\lambda))&=\arg\min(-a_i,\ \lambda+b_i-3a_i,\ -\lambda,\ b_i-2a_i),\\
\label{e:y/v:1}
T(z^{\lambda},v^i(\lambda))&=\arg\min(-b_i,\ \lambda-b_i-a_i,\ -\lambda,\ -b_i),\\
\label{e:y/w:1}
T(z^{\lambda},w^i(\lambda))&=\arg\min(-c_i,\ \lambda-b_i-c_i,\ -\lambda,\ -a_i).
\end{align}
In \eqref{e:y/v:1} and \eqref{e:y/w:1} we used that $2b_i=a_i+c_i$. The inequalities
$a_i\leq\lambda\leq b_i$ imply that
\begin{equation}
\label{ineqs:sort}
-\lambda\leq -a_i\leq b_i-2a_i\leq \lambda+b_i-3a_i,
\end{equation}
hence the minimum in \eqref{e:y/u:1} is attained by the $3$rd component. Analogously, the minimum
in \eqref{e:y/v:1} is attained by the $4$th and $1$st components, 
and the minimum in \eqref{e:y/w:1} is attained
by the $2$nd component. By Theorem \ref{Vorobyov} $z^{\lambda}\in\spann(C(\lambda))$.

\msn {\em Case 4.} $b_i\leq\lambda\leq c_i$. We take
\begin{equation}
\label{ydef2}
\begin{array}{c@{{}\quad{}}ccc}
z^{\lambda}=(0 & c_i & 0 & c_i)^T,
\end{array}
\end{equation}
By Lemma \ref{spannd} $z^{\lambda}\in\spann(D)$, and we
claim again that $z^{\lambda}\in\spann(C(\lambda))$.
We compute
\begin{align}
\label{e:y/v:2}
T(z^{\lambda},v^i(\lambda))&=\arg\min(-b_i,\ c_i-2b_i,\ -\lambda,\ -\lambda),\\
\label{e:y/w:2}
T(z^{\lambda},w^i(\lambda))&=\arg\min(-c_i,\ -c_i,\ -\lambda,\ c_i-b_i-\lambda).
\end{align}
We observe that the minimum in \eqref{e:y/v:2} is attained by 
the $3$rd and $4$th components,
while the minimum in \eqref{e:y/w:2} is attained by the $1$st and $2$nd components. 
The claim follows by Theorem \ref{Vorobyov}.
\end{Proof}

\section*{ACKNOWLEDGEMENT}
\small
This work was supported by EPSRC grant RRAH12809 and 
RFBR grant 08-01-00601.\\
The author is grateful to Peter Butkovi\v{c} for valuable discussions
concerning this paper and the two-sided
eigenproblem, to St\'{e}phane Gaubert for showing him an example that the
spectrum of this eigenproblem can be disconnected and contain both 
an interval and a point, and to Jean-Jacques Loiseau for suggesting
references~\cite{Bru-70} and~\cite{Gan:59}. The author wishes to thank the anonymous referees for 
their careful reading and useful remarks.

\makesubmdate

\makecontacts


\begin{thebibliography}{000}
\bibitem{ABG-06}
M.~Akian, R.~Bapat, and S.~Gaubert:
\newblock Max-plus algebras.
\newblock In L.~Hogben, editor, {\em Handbook of Linear Algebra}, volume~39 of
  {\em Discrete Mathematics and Its Applications}, chapter~25. Chapman and
  Hall/CRC, 2006.

\bibitem{BCOQ}
F.L. Baccelli, G.~Cohen, G.-J. Olsder, and J.-P. Quadrat:
\newblock {\em Synchronization and Linearity: an Algebra for Discrete Event
  Systems}.
\newblock Wiley, 1992.

\bibitem{BV-07}
P.A. Binding and H.~Volkmer:
\newblock A generalized eigenvalue problem in the max algebra.
\newblock {Linear Algebra Appl., {\mi 422} (2007), pp. 360--371.}

\bibitem{Bru-70}
P. Brunovsky:
\newblock A classification of linear controllable systems.
\newblock {Kybernetika, {\mi 6} (1970), pp. 173-188.}

\bibitem{Bur:91}
S.M. Burns:
\newblock Performance analysis and optimization of
asynchronous circuits,
\newblock PhD Thesis, California Institute of Technology, 1991.

\bibitem{But-03}
P.~Butkovi{\v{c}}:
\newblock Max-algebra: the linear algebra of combinatorics?
\newblock {Linear Algebra Appl., {\mi 367} (2003), pp. 313--335.}

\bibitem{But:10}
P.~Butkovi{\v{c}}.
\newblock {\em Max-linear Systems: Theory and Algorithms.}
\newblock Springer, 2010.

\bibitem{Coc-98}
J.~Cochet-Terrasson, G.~Cohen, S.~Gaubert, M.M. Gettrick, and J.P. Quadrat.
\newblock Numerical computation of spectral elements in max-plus algebra.
\newblock In {\em Proceedings of the IFAC conference on systems structure and
  control}, pages 699--706, IRCT, Nantes, France, 1998.


\bibitem{CG:79}
R.~A. Cuninghame-Green.
\newblock {\em Minimax Algebra}, volume 166 of {\em Lecture Notes in Economics
  and Mathematical Systems}.
\newblock Springer, Berlin, 1979.

\bibitem{CGB-03}
R.A. Cuninghame-Green and P.~Butkovi{\v{c}}:
\newblock The equation $A\otimes x=B\otimes y$ over (max,+).
\newblock {Theoretical Computer Science, {\mi 293} (2003), pp. 3--12.}

\bibitem{CGB-08}
R.A. Cuninghame-Green and P.~Butkovi{\v{c}}.
\newblock Generalised eigenproblem in max algebra.
\newblock In {\em Proceedings of the 9th International Workshop WODES 2008},
  pages 236--241, 2008.

\bibitem{ED-01}
L.~Elsner and P.~van~den Driessche:
\newblock Modifying the power method in max algebra.
\newblock {Linear Algebra Appl., {\mi 332-334}, 2001, pp. 3--13.}

\bibitem{Gan:59}
F.R.~Gantmacher:
\newblock {The theory of matrices, Chelsea, 1959.}


\bibitem{GS-10}
S.~Gaubert and S.~Sergeev: The level set method for the two-sided
eigenproblem. E-print \url{http://arxiv.org/pdf/1006.5702}

\bibitem{HOW:05}
B.~Heidergott, G.-J. Olsder, and J.~van~der Woude:
\newblock {\em Max-plus at Work}.
\newblock Princeton Univ. Press, 2005.


\bibitem{McD+98}
J.~J. McDonald, D.~D. Olesky, H.~Schneider, M.~J. Tsatsomeros, and P.~van~den
  Driessche:
\newblock Z-pencils.
\newblock {Electronic J. Linear Algebra, {\mi 4}, 1998, 32--38.}

\bibitem{MNV-08}
V.~Mehrmann, R.~Nabben, and E.~Virnik.
\newblock Generalization of Perron-Frobenius theory to matrix pencils.
\newblock {Linear Algebra Appl., {\mi 428}, 2008, 20--38.}

\end{thebibliography}
\end{document}